\documentclass[runningheads,orivec]{llncs}
\usepackage[T1]{fontenc}
\usepackage{graphicx}
\usepackage{subcaption}
\usepackage{tabularx}
\usepackage[table]{xcolor}
\usepackage{hhline}

\usepackage{amsmath,amssymb}

\begin{document}
\title{Optimization Models for the Quadratic Traveling Salesperson Problem}
%
%
\author{Yuxiao Chen \and
Nivetha Sathish \and
Anubhav Singh \and \\
Ryo Kuroiwa\orcidID{0000-0002-3753-1644} \and
J. Christopher Beck\orcidID{0000-0002-4656-8908}
}
\authorrunning{Chen et al.}
%
\institute{Department of Mechanical and Industrial Engineering, University of Toronto, \\
5 King’s College Road, Toronto, M5S 3G8, Ontario, Canada
\email{\{yuxiao.chen,nivetha.sathish\}@mail.utoronto.ca, anubhav.singh@utoronto.ca, ryo.kuroiwa@mail.utoronto.ca, jcb@mie.utoronto.ca}}
\maketitle              
\begin{abstract}
The quadratic traveling salesperson problem (QTSP) is a generalization of the traveling salesperson problem, in which all triples of consecutive customers in a tour determine the travel cost.
We propose compact optimization models for QTSP in mixed-integer quadratic programming (MIQP), mixed-integer linear programming (MILP), constraint programming (CP), and domain-independent dynamic programming (DIDP).
Our experimental results demonstrate that the DIDP model performs better than other approaches in optimality gap and solution quality when the problem size is large enough.

\end{abstract}
\section{Introduction}

In the Quadratic Traveling Salesperson Problem (QTSP), a set of customers $N = \{ 0, ..., n-1 \}$ is given.
A solution is a tour that visits all customers exactly once and returns to the starting location.
The travel cost is defined by every triple of consecutive customers visited by the tour.
Let $\sigma(i)$ be the $i$-th customer visited, i.e., a tour is represented as $\langle \sigma(0), ..., \sigma(n-1) \rangle$.
Assuming $\sigma(-1) = \sigma(n-1)$ and $\sigma(n) = \sigma(0)$, the cost of the tour is defined as
\begin{equation}
    \sum_{i=0}^{n-1} c_{\sigma(i-1) \sigma(i) \sigma(i+1)}.
\end{equation}

The Quadratic Traveling Salesperson Problem (QTSP) is first defined by Aggrawal et al. \cite{Aggrawal2000} as a variant of the Traveling Salesperson Problem (TSP) that minimizes the angle cost in a tour, the so-called Angle-TSP. The angle cost represents the energy consumption of robots changing directions during the tour, where a larger turning angle results in a higher energy consumption. Thus, the QTSP has been studied for minimizing the energy consumption of robots \cite{Jager2008,Oswin2017,Pham2023,Stanek2019}. The work of Fischer et al. \cite{Fischer2014} is motivated by finding the optimal Permuted Markov model \cite{ellrott2002identifying} or the optimal Permuted
Variable Length Markov model \cite{zhao2005} for a given set of DNA sequences. The authors present a transformation of the bioinformatics problem into a QTSP.


Previous works have proposed specialized algorithms to solve QTSP, including heuristic algorithms \cite{Fischer2014,Pham2023,Stanek2019}, branch-and-bound algorithms \cite{Fischer2014,Jager2008}, and branch-and-cut algorithms \cite{Fischer2014,Jager2008,Oswin2017}. However, none of these works compare the performance of general-purpose solvers for mathematical optimization programs.
While branch-and-cut algorithms are based on integer linear programming (ILP) models and use off-the-shelf solvers as subroutines, they require the implementation of separation algorithms that lazily add constraints to the models and are, hence, specific to QTSP. 

In this paper, we implement and empirically evaluate compact optimization models of QTSP 
in four different paradigms: mixed-integer linear programming (MILP), mixed-integer quadratic programming (MIQP), constraint programming (CP), and domain-independent dynamic programming (DIDP). We solve these models using off-the-shelf solvers directly without any specialized algorithms.

\section{Optimization Models}
In this section, we develop MILP and MIQP models of QTSP, which are built upon existing TSP models. A previous work~\cite{FISCHER201697} proposes similar MILP and MIQP models but uses specialized algorithms to solve them. 
We also present novel CP and dynamic programming formulations of QTSP, which, to our knowledge, have not been investigated previously.


\subsection{Mixed-Integer Linear Programming Model}
Our MILP model is based on a compact MILP model for the traveling salesperson problem (TSP) \cite{DESROCHERS1991}.
We use a binary decision variable $x_{ij}$ to represent visiting customer $j$ after customer $i$, another binary decision variable $y_{ijk}$ to indicate the locations $i, j, k$ are visited consecutively, and an integer decision variable $u_{i}$ to denote the position of customer $i$ in the tour.


\begin{align}
    \min & \sum_{i \in N} \sum_{j \in N \setminus \{ i \}} \sum_{k \in N \setminus \{ i, j \}} c_{ijk} y_{ijk} \label{eqn:milp_obj}\\
    & \sum_{j \in N \setminus \{ i \} } x_{ij} = \sum_{j \in N \setminus \{ i \} } x_{ji} = 1 & \forall i \in N \label{eqn:milp1}\\
    & u_i - u_j + (n-1) x_{ij} + (n-3) x_{ji} \leq n-2 & \begin{array}{ll}
         &  \forall i \in N\setminus \{0\},\\
         &  \forall j \in N\setminus \{0, i\}
    \end{array}  \label{eqn:milp2}\\
    & x_{ij} = \sum_{k \in N \setminus \{ i, j \}} y_{ijk} = \sum_{k \in N \setminus \{ i, j \}} y_{kij} & \forall i \in N, \forall j \in N \setminus \{ i \} \label{eqn:milp3}
\end{align}
\begin{align}
    & x_{ij} \in \{ 0, 1 \} & \forall i \in N, \forall j \in N \setminus \{ i \} \label{eqn:milp4}\\
    & 1 \leq u_{i} \leq n-1 & \forall i \in N \setminus \{0\} \label{eqn:milp5}\\
    & y_{ijk} \in \{ 0, 1 \} & \begin{array}{ll}
         & \forall i \in N, \forall j \in N \setminus \{ i \}, \\
         & \forall k \in N \setminus \{ i, j \}. 
    \end{array}  \label{eqn:milp6}
\end{align}

Constraints~\eqref{eqn:milp1} ensures that each customer has exactly one incoming and outgoing edge.
Constraints~\eqref{eqn:milp2} are the subtour elimination constraints that are introduced by Desrochers and Laporte (DL) \cite{DESROCHERS1991}. They eliminate subtours by assigning a position to each customer in the solution tour, so every customer must have a position index greater than the index of the previous customer by 1, except for customer 0, who has a position index 0. Thus, the selected edges must form a cycle that visits all customers exactly once. We have compared the performance of the MILP models with the DL subtour elimination constraints, Miller-Tucker-Zemlin subtour elimination constraints \cite{MILLER1960}, and the flow-based subtour elimination constraints \cite{Gavish1978}, and the model with DL constraints performs the best. Constraints~\eqref{eqn:milp3} guarantee the $y_{ijk}$ is 1 if and only if both $x_{ij}$ and $x_{jk}$ are 1, which means the customers $i, j, k$ are visited consecutively in the tour. The objective function~\eqref{eqn:milp_obj} is summing over the cost of all three consecutive visits in the tour. In the end, Constraints~\eqref{eqn:milp4}-\eqref{eqn:milp6} bound the domain of each decision variable.

\subsection{Mixed-Integer Quadratic Programming Model}

Our MIQP model is same as the MILP model, except that the MIQP model discards the variables $y_{ijk}$ and uses the variables $x_{ij}$ directly in the objective function, which results in using quadratic terms in the objective function.

\begin{align} 
    \min & \sum_{i \in N} \sum_{j \in N \setminus \{ i \}} \sum_{k \in N \setminus \{ i, j \}} c_{ijk} x_{ij} x_{jk} \label{eqn:miqp_obj}\\
    & \sum_{j \in N \setminus \{ i \} } x_{ij} = \sum_{j \in N \setminus \{ i \} } x_{ji} = 1 & \forall i \in N \label{eqn:miqp1} \\
    & u_i - u_j + (n-1) x_{ij} + (n-3) x_{ji} \leq n-2 & \begin{array}{ll}
         &  \forall i \in N\setminus \{0\},\\
         &  \forall j \in N\setminus \{0, i\}
    \end{array}  \label{eqn:miqp2} \\ 
    & x_{ij} \in \{ 0, 1 \} & \forall i \in N, \forall j \in N \setminus \{ i \} \label{eqn:miqp3}\\
    & 1 \leq u_{i} \leq n-1 & \forall i \in N \setminus \{0\}. \label{eqn:miqp4} 
\end{align}
Constraints~\eqref{eqn:miqp1},~\eqref{eqn:miqp2},~\eqref{eqn:miqp3}, and~\eqref{eqn:miqp4} are the same as Constraints~\eqref{eqn:milp1},~\eqref{eqn:milp2},~\eqref{eqn:milp4}, and~\eqref{eqn:milp5} respectively. The objective function~\eqref{eqn:miqp_obj} includes the cost for consecutively visiting $i,j,k$ if and only if $x_{ij}$ and $x_{jk}$ are both 1.

\subsection{Constraint Programming Model}

In our CP model, we use an integer decision variable $x_i$ to represent the $i$-th customer visited by a tour.
\begin{align}
    \min\ & c_{ x_{n-1} x_0 x_1} + \sum_{i = 0}^{n-3} c_{x_{i} x_{i+1} x_{i+2}} + c_{x_{n-2} x_{n-1} x_0} \label{eqn:cp:objective} \\
    & \mathsf{all\_different}(x_0, ..., x_{n-1}) \label{eqn:cp:all-different} \\
    & x_0 = 0 \label{eqn:cp:depot-location} \\
    & x_i \in N & i = 0, ..., n-1.\label{eqn:cp:domain}
\end{align}
The objective function uses element expressions, which use the values of decision variables as indices of $c$.
Constraint~\eqref{eqn:cp:all-different} ensures that all customers are visited and each customer is visited exactly once. Constraint~\eqref{eqn:cp:depot-location} forces all feasible tours to start at location 0, which reduces the symmetry without loss of generality 
since rooting a customer to a position does not alter the space of possible tour cycles.
Constraints~\eqref{eqn:cp:domain} state the domain of each decision variable.

\subsection{Domain-Independent Dynamic Programming Model}

DIDP is a recently proposed model-based paradigm for dynamic programming (DP) \cite{Kuroiwa2023CAASDy,Kuroiwa2024Journal}.
We can solve a DIDP model formulated in the modeling language, Dynamic Programming Description Language (DyPDL), using general-purpose solvers.
Here, we present a Bellman equation \cite{Bellman1957} for our DIDP model, which can be implemented with DyPDL.

Without loss of generality, we assume that a tour starts from and returns to customer $0$.
In our DIDP model, a problem is represented by four state variables:
$U$ is the set of unvisited customers;
$i$ is the previous customer visited;
$j$ is the current customer;
and $f$ is the first customer visited after $0$.
We consider the next customer to visit from the current state to minimize the total travel cost.
Let $V(U, i, j, f)$ be the optimal cost from the current state.
The DIDP model is represented as follows:
\begin{align}
    \text{compute } & V(N \setminus \{ 0 \}, 0, 0, 0) \label{eqn:dp:objective} \\
    & V(U, i, j, f) = \begin{cases}
        \min_{k \in U} V(U \setminus \{ k \}, 0, k, k) & \text{if } j = 0 \\
        \min_{k \in U} c_{ijk} + V(U \setminus \{ j \}, j, k, f) & \text{if } j \neq 0 \land U \neq \emptyset \\
        c_{j,0,f} + c_{i,j,0} & \text{if } j \neq 0 \land U = \emptyset
    \end{cases} \label{eqn:dp:transitions} \\
    & V(U, i, j, f) \geq \max \left\{ \begin{array}{c}
        \sum_{k \in U \cup \{ f, 0 \}} \min_{l \in N \setminus \{ k \}, m \in N \setminus \{ k, l \} } c_{lmk}, \\ 
        \sum_{k \in U \cup \{ j, 0 \}} \min_{l \in N \setminus \{ k \}, m \in N \setminus \{ k, l \} } c_{lkm}, \\ 
        \sum_{k \in U \cup \{ i, j \}} \min_{l \in N \setminus \{ k \}, m \in N \setminus \{ k, l \} } c_{klm}. \\ 
    \end{array} \right. \label{eqn:dp:dual}
\end{align}

Objective~\eqref{eqn:dp:objective} states that the objective of the model is to compute the optimal cost for the original problem.
Since we assume that $0$ is visited first, the original problem is represented by a state where $U = N \setminus \{ 0 \}$ and $i = 0$.
To represent that the first customer after $0$ is not decided in the original problem, we use $j=f=0$.

Equation~\eqref{eqn:dp:transitions} defines state transitions that transform a state $(U, i, j, f)$ into another state.
When $j = 0$, we decide $f$, the first customer after $0$.
Otherwise, if $U \neq \emptyset$, we visit one customer $k$, and the cost is computed as the sum of the travel cost and the cost of the resulting state.
If all customers are visited, no further state transition is possible, and the cost of a state is the travel cost to visit $f$ and $0$.

Inequality~\eqref{eqn:dp:dual} is a dual bound function, which defines a lower bound on the optimal cost for a state.
We can underestimate the cost to visit a customer $k$ by $\min_{l \in N \setminus \{ k \}, m \in N \setminus \{ k, l \}} c_{lmk}$.
The first line bounds the total cost to visit all unvisited customers using this estimation.
In addition to the set of unvisited customers, it also consider the cost to visit $f$ and $0$.
Similarly, the second line underestimates the cost to visit a customer after $k$, and the third line underestimates the cost to visit a customer when $k$ is previously visited.

\section{Experimental Evaluation}
In our experiments, we used the benchmark instances for QTSP introduced in Stan\v{e}k et al. \cite{Stanek2019}. 
The benchmark contains problem instances with $n$ customers, where $n\in \{5, 10, 15, ..., 200\}$. With each number of customers, there are 10 randomly generated maps, where the locations of the customers are uniformly distributed in the $\{0,...,500\}\times \{0,...,500\}$ grid. The benchmark has two problem instances for each map: an AngleTSP instance and an AngleDistanceTSP instance, where the difference is on the cost functions.
\begin{itemize}
    \item \textbf{AngleTSP instances} use the turning angle as the cost. Specifically, suppose the vehicle visits location $i, j, \text{and } k$ in order, then the cost is the turning angle $\alpha_{ijk}$ between the vectors $\vec{ij}$ and $\vec{jk}$, multiplied by 1000, and rounded to 12 decimal places.
    \item \textbf{AngleDistanceTSP instances} have a cost function that combines the turning angle with the Euclidean distances between the points in a weighted sum. Let $d_{ij}$ represent the Euclidean distance between $i$ and $j$, and $\rho\in \mathbb{R}^+_0$ be a weighting parameter, then the cost of visiting locations $i, j, k$ in order is.
    \[c_{ijk} = 100\left(\rho\cdot \alpha_{ijk} + \frac{d_{ij}+d_{jk}}{2}\right),\]
    where $\alpha_{ijk}$ is the turning angle used in the AngleTSP instances. Notice as $\rho \rightarrow \infty$, the instance is similar to AngleTSP instances, and as $\rho \rightarrow 0$, the instance is similar to the standard TSP instances. In this benchmark, $\rho$ is set to 40 for all instances.
    
\end{itemize}
    
We used four different models and solvers to solve the benchmark instances: the MIQP model and MILP model with Gurobi 11.0.3 \cite{gurobi}, the DIDP model with the Complete Anytime Beam Search (CABS) \cite{Zhang1998,Kuroiwa2023Anytime} in DIDPPy 0.8.0, and the CP model with IBM ILOG CP Optimizer 22.1.1 \cite{Laborie2018}. All experiments are given a time limit of 1800 seconds and 8GB memory, and performing on an Intel Xeon Gold 6148 core at 2.4GHz using GNU Parallel \cite{Tange2011}. 

\begin{table}[htbp]
    \fontsize{8}{10}\selectfont
    \centering
    \begin{tabularx}{\textwidth}{c || c | c | c | c | c | c | c | c | c | c | c | c | c | c | c | c | c | c | c | c}
        \hline
         & \multicolumn{20}{c}{Number of Locations} \\
        \hline
        Solver & 5 & 10 & 15 & 20 & 25 & 30 & 35 & 40 & 45 & 50 & 55 & 60 & 65 & 70 & 75 & 80 & 85 & 90 & 95 & 100 \\ 
        \hhline{=||=|=|=|=|=|=|=|=|=|=|=|=|=|=|=|=|=|=|=|=}
        DIDP & 0 & 0 & 0 & \textcolor{red}{9} & \textcolor{red}{10} & 0 & 0 & \textcolor{red}{3} & 0 & 0 & 0 & 0 & 0 & 0 & 0 & 0 & 0 & 0 & 0 & 0 \\ 
        \hline
        MILP & 0 & 0 & 0 & 0 & 0 & 0 & 0 & 0 & 0 & 0 & 0 & 0 & 0 & 0 & 0 & 0 & 0 & 0 & 0 & 0 \\
        \hline
        MIQP & 0 & 0 & 0 & 0 & 0 & 0 & 0 & 0 & 0 & 0 & 0 & 0 & 0 & 0 & 0 & 0 & 0 & 0 & 0 & \textcolor{red}{10} \\
        \hline
        CP & 0 & 0 & 0 & 0 & 0 & 0 & 0 & 0 & 0 & 0 & 0 & 0 & 0 & 0 & 0 & 0 & 0 & 0 & 0 & 0 \\
        \hline \hline
         & \multicolumn{20}{c}{Number of Locations} \\
         \hline
        Solver & 105 & 110 & 115 & 120 & 125 & 130 & 135 & 140 & 145 & 150 & 155 & 160 & 165 & 170 & 175 & 180 & 185 & 190 & 195 & 200 \\
        \hhline{=||=|=|=|=|=|=|=|=|=|=|=|=|=|=|=|=|=|=|=|=}
        DIDP & 0 & 0 & 0 & 0 & 0 & 0 & 0 & 0 & 0 & 0 & 0 & 0 & 0 & 0 & 0 & 0 & 0 & 0 & 0 & 0 \\
        \hline
        MILP & 0 & 0 & 0 & 0 & 0 & 0 & 0 & 0 & 0 & 0 & 0 & 0 & 0 & 0 & \textcolor{red}{10} & \textcolor{red}{10} & \textcolor{red}{10} & \textcolor{red}{10} & \textcolor{red}{10} & \textcolor{red}{10} \\
        \hline
        MIQP & \textcolor{red}{10} & 0 & 0 & 0 & 0 & 0 & 0 & 0 & 0 & 0 & 0 & 0 & 0 & 0 & 0 & 0 & 0 & 0 & 0 & \textcolor{red}{10} \\
        \hline
        CP & 0 & \textcolor{red}{10} & \textcolor{red}{10} & \textcolor{red}{10} & \textcolor{red}{10} & \textcolor{red}{10} & \textcolor{red}{10} & \textcolor{red}{10} & \textcolor{red}{10} & \textcolor{red}{10} & \textcolor{red}{10} & \textcolor{red}{10} & \textcolor{red}{10} & \textcolor{red}{10} & \textcolor{red}{10} & \textcolor{red}{10} & \textcolor{red}{10} & \textcolor{red}{10} & \textcolor{red}{10} & \textcolor{red}{10} \\ 
        \hline
    \end{tabularx}
    \caption{Number of AngleTSP instances that each solver has memory-out} issue.\label{table:memout_angle}
\end{table}

\begin{table}[htbp]
    \fontsize{8}{10}\selectfont
    \centering
    \begin{tabularx}{\textwidth}{c || c | c | c | c | c | c | c | c | c | c | c | c | c | c | c | c | c | c | c | c}
        \hline
         & \multicolumn{20}{c}{Number of Locations} \\
        \hline
        Solver & 5 & 10 & 15 & 20 & 25 & 30 & 35 & 40 & 45 & 50 & 55 & 60 & 65 & 70 & 75 & 80 & 85 & 90 & 95 & 100 \\ 
        \hhline{=||=|=|=|=|=|=|=|=|=|=|=|=|=|=|=|=|=|=|=|=}
        DIDP & 0 & 0 & 0 & 0 & 0 & 0 & 0 & 0 & 0 & 0 & 0 & 0 & 0 & 0 & 0 & 0 & 0 & 0 & 0 & 0 \\ 
        \hline
        MILP & 0 & 0 & 0 & 0 & 0 & 0 & 0 & 0 & 0 & 0 & 0 & 0 & 0 & 0 & 0 & 0 & 0 & 0 & 0 & 0 \\
        \hline
        MIQP & 0 & 0 & 0 & 0 & 0 & 0 & 0 & 0 & 0 & 0 & 0 & 0 & 0 & 0 & 0 & 0 & 0 & 0 & 0 & \textcolor{red}{8} \\
        \hline
        CP & 0 & 0 & 0 & 0 & 0 & 0 & 0 & 0 & 0 & 0 & 0 & 0 & 0 & 0 & 0 & 0 & 0 & 0 & 0 & 0 \\
        \hline \hline
         & \multicolumn{20}{c}{Number of Locations} \\
         \hline
        Solver & 105 & 110 & 115 & 120 & 125 & 130 & 135 & 140 & 145 & 150 & 155 & 160 & 165 & 170 & 175 & 180 & 185 & 190 & 195 & 200 \\
        \hhline{=||=|=|=|=|=|=|=|=|=|=|=|=|=|=|=|=|=|=|=|=}
        DIDP & 0 & 0 & 0 & 0 & 0 & 0 & 0 & 0 & 0 & 0 & 0 & 0 & 0 & 0 & 0 & 0 & 0 & 0 & 0 & 0 \\
        \hline
        MILP & 0 & 0 & 0 & 0 & 0 & 0 & 0 & 0 & 0 & 0 & 0 & 0 & 0 & 0 & \textcolor{red}{10} & \textcolor{red}{10} & \textcolor{red}{10} & \textcolor{red}{10} & \textcolor{red}{10} & \textcolor{red}{10} \\
        \hline
        MIQP & \textcolor{red}{5} & 0 & 0 & 0 & 0 & 0 & 0 & 0 & 0 & 0 & 0 & 0 & 0 & 0 & 0 & 0 & 1 & 0 & 0 & \textcolor{red}{10} \\
        \hline
        CP & 0 & \textcolor{red}{10} & \textcolor{red}{10} & \textcolor{red}{10} & \textcolor{red}{10} & \textcolor{red}{10} & \textcolor{red}{10} & \textcolor{red}{10} & \textcolor{red}{10} & \textcolor{red}{10} & \textcolor{red}{10} & \textcolor{red}{10} & \textcolor{red}{10} & \textcolor{red}{10} & \textcolor{red}{10} & \textcolor{red}{10} & \textcolor{red}{10} & \textcolor{red}{10} & \textcolor{red}{10} & \textcolor{red}{10} \\ 
        \hline
    \end{tabularx}
    \caption{Number of AngleDistanceTSP instances that each solver has memory-out issue.}\label{table:memout_angleDist}
\end{table}

\begin{table}[htbp]
    \fontsize{8}{10}\selectfont
    \centering
    \begin{tabularx}{\textwidth}{c || c | c | c | c | c | c | c | c | c | c | c | c | c | c | c | c | c | c | c | c}
        \hline
         & \multicolumn{20}{c}{Number of Locations} \\
        \hline
        Solver & 5 & 10 & 15 & 20 & 25 & 30 & 35 & 40 & 45 & 50 & 55 & 60 & 65 & 70 & 75 & 80 & 85 & 90 & 95 & 100 \\ 
        \hhline{=||=|=|=|=|=|=|=|=|=|=|=|=|=|=|=|=|=|=|=|=}
        DIDP & \textcolor{red}{10} & \textcolor{red}{10} & \textcolor{red}{10} & 0 & 0 & 0 & 0 & 3 & 0 & 0 & 0 & 0 & 0 & 0 & 0 & 0 & 0 & 0 & 0 & 0 \\ 
        \hline
        MILP & \textcolor{red}{10} & \textcolor{red}{10} & \textcolor{red}{10} & \textcolor{red}{10} & \textcolor{red}{10} & \textcolor{red}{10} & \textcolor{red}{10} & \textcolor{red}{7} & \textcolor{red}{3} & \textcolor{red}{1} & 0 & 0 & 0 & 0 & 0 & 0 & 0 & 0 & 0 & 0 \\
        \hline
        MIQP & \textcolor{red}{10} & \textcolor{red}{10} & 0 & 0 & 0 & 0 & 0 & 0 & 0 & 0 & 0 & 0 & 0 & 0 & 0 & 0 & 0 & 0 & 0 & 0 \\
        \hline
        CP & \textcolor{red}{10} & \textcolor{red}{10} & 0 & 0 & 0 & 0 & 0 & 0 & 0 & 0 & 0 & 0 & 0 & 0 & 0 & 0 & 0 & 0 & 0 & 0 \\
        \hline \hline
         & \multicolumn{20}{c}{Number of Locations} \\
         \hline
        Solver & 105 & 110 & 115 & 120 & 125 & 130 & 135 & 140 & 145 & 150 & 155 & 160 & 165 & 170 & 175 & 180 & 185 & 190 & 195 & 200 \\
        \hhline{=||=|=|=|=|=|=|=|=|=|=|=|=|=|=|=|=|=|=|=|=}
        DIDP & 0 & 0 & 0 & 0 & 0 & 0 & 0 & 0 & 0 & 0 & 0 & 0 & 0 & 0 & 0 & 0 & 0 & 0 & 0 & 0 \\
        \hline
        MILP & 0 & 0 & 0 & 0 & 0 & 0 & 0 & 0 & 0 & 0 & 0 & 0 & 0 & 0 & 0 & 0 & 0 & 0 & 0 & 0 \\
        \hline
        MIQP & 0 & 0 & 0 & 0 & 0 & 0 & 0 & 0 & 0 & 0 & 0 & 0 & 0 & 0 & 0 & 0 & 0 & 0 & 0 & 0 \\
        \hline
        CP & 0 & 0 & 0 & 0 & 0 & 0 & 0 & 0 & 0 & 0 & 0 & 0 & 0 & 0 & 0 & 0 & 0 & 0 & 0 & 0 \\ 
        \hline
    \end{tabularx}
    \caption{Number of AngleTSP instances that each solver has found the optimal solution.}\label{table:optimal_angle}
\end{table}

\begin{table}[htbp]
    \fontsize{8}{10}\selectfont
    \centering
    \begin{tabularx}{\textwidth}{c || c | c | c | c | c | c | c | c | c | c | c | c | c | c | c | c | c | c | c | c}
        \hline
         & \multicolumn{20}{c}{Number of Locations} \\
        \hline
        Solver & 5 & 10 & 15 & 20 & 25 & 30 & 35 & 40 & 45 & 50 & 55 & 60 & 65 & 70 & 75 & 80 & 85 & 90 & 95 & 100 \\ 
        \hhline{=||=|=|=|=|=|=|=|=|=|=|=|=|=|=|=|=|=|=|=|=}
        DIDP & \textcolor{red}{10} & \textcolor{red}{10} & \textcolor{red}{10} & \textcolor{red}{10} & \textcolor{red}{10} & \textcolor{red}{6} & 0 & 0 & 0 & 0 & 0 & 0 & 0 & 0 & 0 & 0 & 0 & 0 & 0 & 0 \\ 
        \hline
        MILP & \textcolor{red}{10} & \textcolor{red}{10} & \textcolor{red}{10} & \textcolor{red}{10} & \textcolor{red}{10} & \textcolor{red}{10} & \textcolor{red}{10} & \textcolor{red}{10} & \textcolor{red}{10} & \textcolor{red}{10} & \textcolor{red}{9} & \textcolor{red}{10} & \textcolor{red}{7} & \textcolor{red}{5} & \textcolor{red}{3} & 0 & 0 & 0 & 0 & 0\\
        \hline
        MIQP & \textcolor{red}{10} & \textcolor{red}{10} & \textcolor{red}{2} & 0 & 0 & 0 & 0 & 0 & 0 & 0 & 0 & 0 & 0 & 0 & 0 & 0 & 0 & 0 & 0 & 0 \\
        \hline
        CP & \textcolor{red}{10} & \textcolor{red}{10} & 0 & 0 & 0 & 0 & 0 & 0 & 0 & 0 & 0 & 0 & 0 & 0 & 0 & 0 & 0 & 0 & 0 & 0 \\
        \hline \hline
         & \multicolumn{20}{c}{Number of Locations} \\
         \hline
        Solver & 105 & 110 & 115 & 120 & 125 & 130 & 135 & 140 & 145 & 150 & 155 & 160 & 165 & 170 & 175 & 180 & 185 & 190 & 195 & 200 \\
        \hhline{=||=|=|=|=|=|=|=|=|=|=|=|=|=|=|=|=|=|=|=|=}
        DIDP & 0 & 0 & 0 & 0 & 0 & 0 & 0 & 0 & 0 & 0 & 0 & 0 & 0 & 0 & 0 & 0 & 0 & 0 & 0 & 0 \\
        \hline
        MILP & 0 & 0 & 0 & 0 & 0 & 0 & 0 & 0 & 0 & 0 & 0 & 0 & 0 & 0 & 0 & 0 & 0 & 0 & 0 & 0 \\
        \hline
        MIQP & 0 & 0 & 0 & 0 & 0 & 0 & 0 & 0 & 0 & 0 & 0 & 0 & 0 & 0 & 0 & 0 & 0 & 0 & 0 & 0 \\
        \hline
        CP & 0 & 0 & 0 & 0 & 0 & 0 & 0 & 0 & 0 & 0 & 0 & 0 & 0 & 0 & 0 & 0 & 0 & 0 & 0 & 0 \\ 
        \hline
    \end{tabularx}
    \caption{Number of AngleDistanceTSP instances that each solver has found the optimal solution.}\label{table:optimal_angleDist}
\end{table}

Tables \ref{table:memout_angle} and \ref{table:memout_angleDist} show the number of AngleTSP and AngleDistance instances that each solver cannot solve within the 8GB memory limit, and the solvers with the most memory-out instances are marked red. We observe that all solvers except the DIDP solver experience a memory-out issue when the problem size is large enough. However, DIDP solver reaches the memory limit while solving some small AngleTSP instances (20 to 25 customers), which is because that the CABS algorithm expands states faster when the problem size is smaller, so the algorithm consumes more memory within the time limit compared to the larger instances. We also notice that the Gurobi solver reaches the memory limit when solving both the AngleTSP and AngleDistanceTSP instances with 100 and 105 customers using the MIQP model. However, we are unable to explain this from the detailed logs of the solver runs.

Tables \ref{table:optimal_angle} and \ref{table:optimal_angleDist} show the number of AngleTSP and AngleDistanceTSP instances that each solver can find and prove the optimal solution, and the solvers that have proved the most optimal solutions are marked in red. We observe the Gurobi solver with MILP model performs the best, and the DIDP model performs better than MIQP and CP models.

\begin{figure}
\centering
\begin{subfigure}{.5\textwidth}
  \centering
  \includegraphics[width=\linewidth]{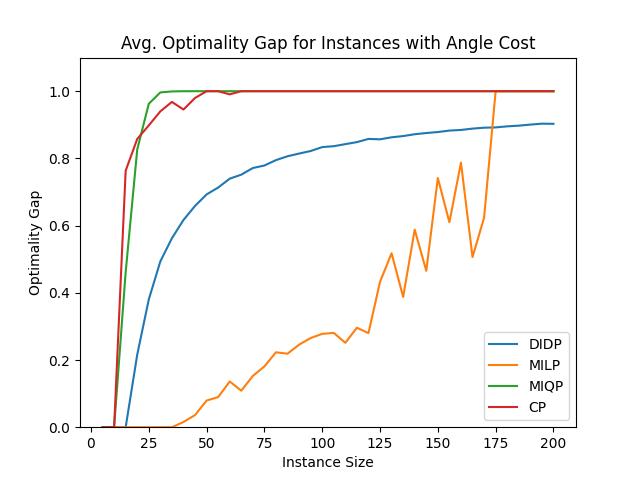}
  \captionsetup{width=.9\linewidth}
  \caption{The average \textit{optimality gap} found by each solver for the AngleTSP instances.}
  \label{fig:opt_gap_angle}
\end{subfigure}%
\begin{subfigure}{.5\textwidth}
  \centering
  \includegraphics[width=\linewidth]{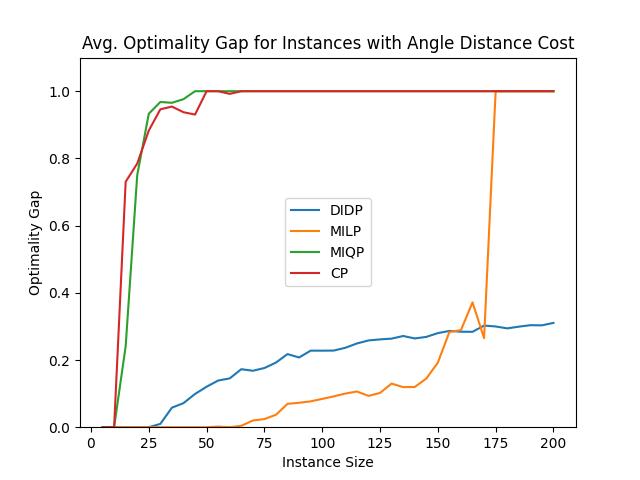}
  \captionsetup{width=.9\linewidth}
  \caption{The average \textit{optimality gap} found by each solver for the AngleDistance-TSP instances.}
  \label{fig:opt_gap_angle_dist}
\end{subfigure}
\caption{The plots of average \textit{optimality gap} found by each solver.}
\label{fig:opt_gap}
\end{figure}

Figure \ref{fig:opt_gap} shows the average \textit{optimality gap} over the 10 instances with the same instance size (number of customers) and cost type (Angle or AngleDistance). The \textit{optimality gap} is calculated by
\[\frac{|\text{Primal Bound} - \text{Dual Bound}|}{\text{Primal Bound}},\]
where the \textbf{Primal Bound} is the cost of the best found solution, which is an upper bound on the optimal solution in a minimization problem, and the \textbf{Dual Bound} is the best lower bound on the optimal solution proved by the solver. Notice 0 is always a lower bound on the cost of the optimal solution in a QTSP, so the maximum \textit{optimality gap} is 1. Thus, the \textit{optimality gap} of the instances that have no feasible solution found is set to 1.
From Figure \ref{fig:opt_gap}, we observe that the MILP model has the smallest \textit{optimality gap} when solving small instances, which means it is closer to proving the optimality of a solution, compared to the other models. However, when the problem size is large enough ($\geq175$ customers), the MILP model cannot find any primal or dual bound within the memory limit, and the DIDP model has the best \textit{optimality gap} for these large instances. The same behavior is observed for solving both the AngleTSP instances (Figure \ref{fig:opt_gap_angle}) and AngleDistanceTSP instances (Figure \ref{fig:opt_gap_angle_dist}). 

\begin{figure}
\centering
\begin{subfigure}[t]{.5\textwidth}
  \centering
  \includegraphics[width=\linewidth]{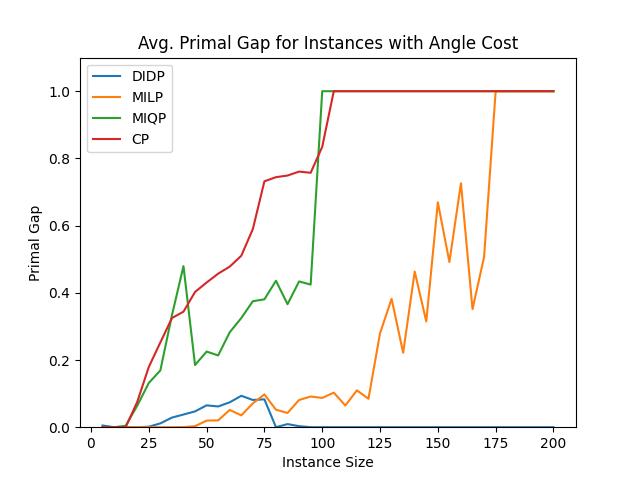}
  \captionsetup{width=.9\linewidth}
  \caption{The average \textit{primal gap} found by each solver for the AngleTSP instances.}
  \label{fig:primal_gap_angle}
\end{subfigure}%
\begin{subfigure}[t]{.5\textwidth}
  \centering
  \includegraphics[width=\linewidth]{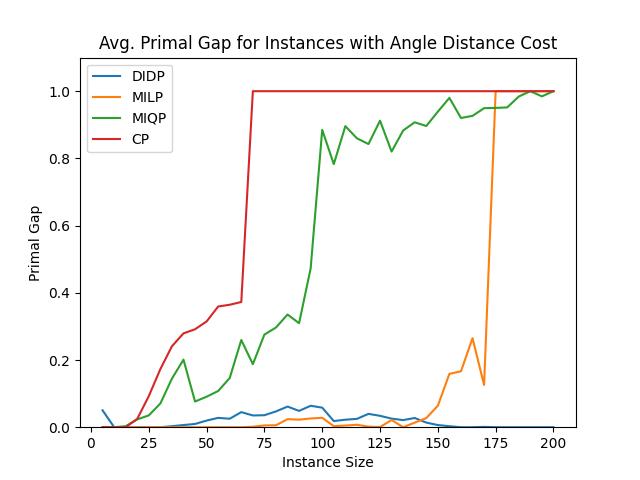}
  \captionsetup{width=.9\linewidth}
  \caption{The average \textit{primal gap} found by each solver for the AngleDistanceTSP instances.}
  \label{fig:primal_gap_angle_dist}
\end{subfigure}
\caption{The plots of average \textit{primal gap} found by each solver.}
\label{fig:primal_gap}
\end{figure}

\begin{figure}
\centering
\begin{subfigure}[t]{.5\textwidth}
  \centering
  \includegraphics[width=\linewidth]{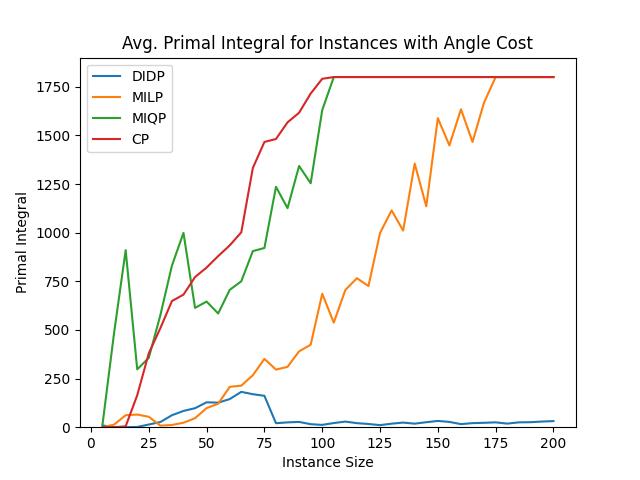}
  \captionsetup{width=.9\linewidth}
  \caption{The average \textit{primal integral} of the search results from each solver for the AngleTSP instances.}
  \label{fig:primal_int_angle}
\end{subfigure}%
\begin{subfigure}[t]{.5\textwidth}
  \centering
  \includegraphics[width=\linewidth]{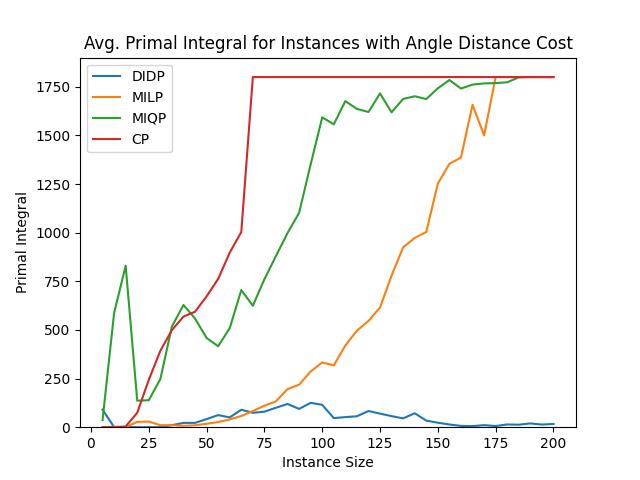}
  \captionsetup{width=.9\linewidth}
  \caption{The average \textit{primal integral} of the search results from each solver for the AngleDistanceTSP instances.}
  \label{fig:primal_int_angle_dist}
\end{subfigure}
\caption{The plots of average \textit{primal integral} produced by each solver.}
\label{fig:primal_int}
\end{figure}

Figure \ref{fig:primal_gap_angle} shows the average \textit{primal gap} of the 10 AngleTSP instances found by each solver, and Figure \ref{fig:primal_gap_angle_dist} shows the average \textit{primal gap} of the AngleDistanceTSP instances. The \textit{primal gap} is calculated by
\[\frac{|\text{Primal Bound} - \text{Best Known Solution}|}{\text{Primal Bound}},\]
where the \textbf{Best Known Solution} is the optimal solution for the instances with $\leq 75$ customers given in the benchmark, and best found solution in all 4 solvers is the \textbf{Best Known Solution} for the large instances. \textit{Primal gap} is a measurement for the solution quality, a smaller \textit{Primal gap} implies the cost of the feasible solution found is smaller. The maximum value of the \textit{primal gap} is 1, when there has no feasible solution found, which means the primal bound has a value of $\infty$. From the Figure \ref{fig:primal_gap}, we observe that the MILP model finds the best feasible solutions for the instances with small number of customers ($\leq 75$ for AngleTSP instances and $\leq 140$ for AngleDistanceTSP instances). As the problem size increases, DIDP model finds the best feasible solution compared to all other solvers. The reason is that the DIDP model for QTSP problem has a feasible solution with any permutation of the visits, and the CABS algorithm performs like a depth-first search when the beam size is small, so it always finds some feasible solution in a short time.
Moreover, Figure \ref{fig:primal_int} shows the average \textit{primal integral} of the search results produced by each solver. \textit{Primal integral} is a measurement for the solution quality and the time they have been found \cite{BERTHOLD2013611}, it is the area under the \textit{primal gap} vs time line, so a smaller \textit{primal integral} indicates a solver finds a better feasible solution faster. From Figure \ref{fig:primal_int}, we observe the DIDP model starts to have the best \textit{primal integral} for the problems with more than 60 customers, which means the DIDP model has the best \textit{primal integral} in more instances compared to the ones that the model has the best \textit{primal gap}. Thus, the DIDP model has an advantage on finding good feasible solutions for the large instances, and also the feasible solutions are found faster compared to the other solvers. 

\section{Conclusion}


The four proposed QTSP models show unique performance profiles on the benchmark instances. The DIDP model seems to scale very well compared to others; all models except DIDP exceed the memory or time limits for large problems without finding any feasible solution. In contrast, the DIDP model always finds a feasible solution quickly, so it outperforms all other approaches in terms of the \textit{optimality gap} and solution quality on large problems. However, DIDP proves optimality only for the three smallest instance sizes. MILP, on the other hand, has the highest count of instances solved to optimality but fails to find a feasible solution for the large instances.
This points to a need to design and develop efficient methods to improve lower bound in DIDP in the future, such as adding more dual bounds.

%
%
%
 \bibliographystyle{splncs04}
 \bibliography{references}

\begin{thebibliography}{10}
\providecommand{\url}[1]{\texttt{#1}}
\providecommand{\urlprefix}{URL }
\providecommand{\doi}[1]{https://doi.org/#1}

\bibitem{Aggrawal2000}
Aggarwal, A., Coppersmith, D., Khanna, S., Motwani, R., Schieber, B.: The angular-metric traveling salesman problem. SIAM Journal on Computing  \textbf{29}(3),  697--711 (2000). \doi{10.1137/S0097539796312721}

\bibitem{Bellman1957}
Bellman, R.: Dynamic Programming. Princeton University Press (1957)

\bibitem{BERTHOLD2013611}
Berthold, T.: Measuring the impact of primal heuristics. Operations Research Letters  \textbf{41}(6),  611--614 (2013). \doi{10.1016/j.orl.2013.08.007}

\bibitem{DESROCHERS1991}
Desrochers, M., Laporte, G.: Improvements and extensions to the {M}iller-{T}ucker-{Z}emlin subtour elimination constraints. Operations Research Letters  \textbf{10}(1),  27--36 (1991). \doi{10.1016/0167-6377(91)90083-2}

\bibitem{ellrott2002identifying}
Ellrott, K., Yang, C., Sladek, F.M., Jiang, T.: Identifying transcription factor binding sites through markov chain optimization. Bioinformatics  \textbf{18}(suppl\_2),  S100--S109 (2002). \doi{10.1093/bioinformatics/18.suppl_2.s100}

\bibitem{Fischer2014}
Fischer, A., Fischer, F., J{\"a}ger, G., Keilwagen, J., Molitor, P., Grosse, I.: Exact algorithms and heuristics for the quadratic traveling salesman problem with an application in bioinformatics. Discrete Applied Mathematics  \textbf{166},  97--114 (2014). \doi{10.1016/j.dam.2013.09.011}

\bibitem{FISCHER201697}
Fischer, A., Fabian~Meier, J., Pferschy, U., Staněk, R.: Linear models and computational experiments for the quadratic {TSP}. Electronic Notes in Discrete Mathematics  \textbf{55},  97--100 (2016). \doi{10.1016/j.endm.2016.10.025}

\bibitem{Gavish1978}
Gavish, B., Graves, S.C.: The travelling salesman problem and related problems. Tech. rep., Operations Research Center, Massachusetts Institute of Technology (1978), {W}orking Paper OR 078-78

\bibitem{gurobi}
{Gurobi Optimization, LLC}: {Gurobi Optimizer Reference Manual} (2023), \url{https://www.gurobi.com}

\bibitem{Jager2008}
J{\"a}ger, G., Molitor, P.: Algorithms and experimental study for the traveling salesman problem of second order. In: Combinatorial Optimization and Applications. pp. 211--224. Springer, Berlin, Heidelberg (2008). \doi{10.1007/978-3-540-85097-7_20}

\bibitem{Kuroiwa2023CAASDy}
Kuroiwa, R., Beck, J.C.: Domain-independent dynamic programming: Generic state space search for combinatorial optimization. In: Proceedings of the 33rd International Conference on Automated Planning and Scheduling (ICAPS). pp. 236--244. {AAAI} Press, Palo Alto, California USA (2023). \doi{10.1609/icaps.v33i1.27200}

\bibitem{Kuroiwa2023Anytime}
Kuroiwa, R., Beck, J.C.: Solving domain-independent dynamic programming problems with anytime heuristic search. In: Proceedings of the 33rd International Conference on Automated Planning and Scheduling (ICAPS). pp. 245--253. {AAAI} Press, Palo Alto, California USA (2023). \doi{10.1609/icaps.v33i1.27201}

\bibitem{Kuroiwa2024Journal}
Kuroiwa, R., Beck, J.C.: Domain-independent dynamic programming. arXiv:2401.13883 [cs.AI] (2024). \doi{10.48550/arXiv.2401.13883}

\bibitem{Laborie2018}
Laborie, P., Rogerie, J., Shaw, P., Vilím, P.: {IBM} {ILOG} {CP} optimizer for scheduling. Constraints  \textbf{23}(2),  210--250 (2018). \doi{10.1007/s10601-018-9281-x}

\bibitem{MILLER1960}
Miller, C.E., Tucker, A.W., Zemlin, R.A.: Integer programming formulation of traveling salesman problems. J. ACM  \textbf{7}(4),  326–329 (oct 1960). \doi{10.1145/321043.321046}, \url{https://doi.org/10.1145/321043.321046}

\bibitem{Oswin2017}
Oswin, A., Fischer, A., Fischer, F., Meier, J.F., Pferschy, U., Pilz, A., Stan\v{e}k, R.: Minimization and maximization versions of the quadratic travelling salesman problem. Optimization  \textbf{66}(4),  521--546 (2017). \doi{10.1080/02331934.2016.1276905}

\bibitem{Pham2023}
Pham, Q.A., Lau, H.C., {H\`a}, M.H., Vu, L.: An efficient hybrid genetic algorithm for the quadratic traveling salesman problem. In: Proceedings of the 33rd International Conference on Automated Planning and Scheduling (ICAPS). pp. 343--351 (2023). \doi{10.1609/icaps.v33i1.27212}

\bibitem{Stanek2019}
Stan\v{e}k, R., Greistorfer, P., Ladner, K., Pferschy, U.: Geometric and {LP}-based heuristics for angular travelling salesman problems in the plane. Computers \& Operations Research  \textbf{108},  97--111 (2019). \doi{10.1016/j.cor.2019.01.016}

\bibitem{Tange2011}
Tange, O.: {GNU} parallel - the command-line power tool. The USENIX Magazine  \textbf{36},  42--47 (2011)

\bibitem{Zhang1998}
Zhang, W.: Complete anytime beam search. In: Proceedings of the 15th National Conference on Artificial Intelligence (AAAI). pp. 425--430. {AAAI} Press (1998)

\bibitem{zhao2005}
Zhao, X., Huang, H., Speed, T.P.: Finding short dna motifs using permuted markov models. Journal of Computational Biology  \textbf{12}(6),  894--906 (2005). \doi{10.1089/cmb.2005.12.894}, \url{https://doi.org/10.1089/cmb.2005.12.894}, pMID: 16108724

\end{thebibliography}
\end{document}